\documentclass[11pt]{article}
\usepackage{amssymb}
\usepackage{latexsym}
\newcommand{\reals}{{\mathbb R}}

\newcommand{\calb}{{\cal B}}

\newcommand{\cale}{{\cal E}}
\newcommand{\calg}{{\cal G}}

\newcommand{\calj}{{\cal J}}
\newcommand{\call}{{\cal L}}

\newcommand{\cals}{{\cal S}}

\newcommand{\del}{\partial}

\newcommand{\frakg}{\mathfrak{g}}
\newcommand{\frakh}{\mathfrak{h}}

\newcommand{\arrows}{\,\lower1pt\hbox{$\longrightarrow$}\hskip-.24in\raise2pt
             \hbox{$\longrightarrow$}\,}
\title{{\bf The Geometry of Momentum}}
\author{Alan
Weinstein\thanks{Research partially supported by NSF Grant
DMS-99-71505.
\newline \mbox{~~~~}MSC2000 Subject Classification Numbers: 53D20
(Primary), 20L05, 53D17, 70H33 (Secondary).
\newline \mbox{~~~~}Keywords: moment map, momentum map,
hamiltonian, symplectic, symmetry, conservation law, groupoid}
\\Department of Mathematics\\ University of California\\ Berkeley, CA
94720 USA\\ {\small(alanw@math.berkeley.edu)}}
\date{}
\begin{document}
\maketitle
\begin{center}
To Jerry Marsden for his 60th birthday
\end{center}
\begin{abstract}
Although the idea of the momentum map associated with a symplectic action of a
group is already contained in work of Lie, the geometry of momentum
maps was not studied extensively until the 1960's. Centering around
the relation between symmetries and conserved quantities, the study of 
momentum maps was very much alive at the end of the 20th century and
continues to this day, with the creation of new notions of symmetry.
A uniform framework for all these momentum maps is still to be found;
groupoids should play an important role in such a framework.

\end{abstract}

\section{Introduction}
\label{sec-intro}
The term {\em momentum} (in Italian, {\em momento})
was introduced by Galileo Galilei as the ``virtue'' of a moving
object which keeps it moving.  He expressed it as the product of
weight and velocity.   Galileo's notion of momentum
was an outgrowth of the medieval {\em impetus} theory of William of
Occam, Jean Buridan, Nicole Oresme, and others (see, for example,
Dugas \cite{du:mecanique}), which deviated from
Aristotle's view that something external to a moving (terrestrial)
body was necessary to keep it in motion.  (Opposition to Aristotle's
view had also been expressed as early as the 5th century by 
John (Philopon) of Alexandria.)

By now, ``momentum'' has a much wider meaning.  Many applications of
the term refer to quantities whose conservation under the time evolution of
a physical system is related to some symmetry of the system.  The
theme of this paper (and the corresponding talk at the conference) is
this relation between momentum and symmetry, especially as the meaning
of the later was extended in the late 20th century.  Though much of
this extension was motivated by quantum mechanics, the emphasis here
will be on momentum in classical mechanics.  Another very important
aspect of our subject which we will not treat in this paper is {\em
symplectic reduction}, which is the simultaneous use of symmetries and
conserved quantities to reduce the dimensionality of a hamiltonian
system.  We refer to \cite{ma-we:comments} for an extensive history of
reduction.  

\section{Momentum in the calculus of variations}
Despite the 
original conception of momentum as {\em proportional} to velocity, we now see
momentum for the most part as {\em dual} to velocity.  This idea
appears clearly in the calculus of variations where, given a
lagrangian $L(t,x,\dot{x})$ for a function $x=(x^1,...,x^n)$ of $t$,
we define the conjugate momenta by $p_i=\del L/\del {\dot x}^i$ and then the
hamiltonian by $H=\sum p_i {\dot x}^i - L$.  One sees from the
way in which the conjugate momenta transform under coordinate changes
(which results in the invariance of $H$ under such changes) that the
conjugate momenta are the components of a vector $p$ which is dual to
the velocity vector ${\dot x}.$  

Of course, if $L=\frac{1}{2}m\sum_i ({\dot x}^i)^2 - V(x),$ then
$p_i=m{\dot x}^i$.  The conflict between lower and upper indices here
results from the use of orthonormal coordinates.  More
invariantly, we should write the first (kinetic energy) part of $L$ as
\begin{equation}
\label{eq-simple}
\frac{1}{2}\sum_{i,j} m_{ij}{\dot x}^i{\dot x}^j,
\end{equation}
 where
$m_{ij} = m \delta_{ij}$.  The ``mass tensor'' $m_{ij}$ may then replaced
in (\ref{eq-simple}) by 
an arbitrary riemannian metric $m_{ij}(x)$ to describe in curvilinear
coordinates what Smale \cite{sm:topologyI} calls a {\em simple
mechanical system}.  For such a system, the relation between momentum
and velocity becomes 
\begin{equation}
\label{eq-psimple}
p_i = \sum_j m_{ij}{\dot x}^j,
\end{equation}
 the mass tensor
now being interpreted as a mapping from velocity (tangent) to momentum
(cotangent) vectors.

If we add the vector potential of a magnetic field to obtain
$$L=\frac{1}{2}\sum_{i,j} m_{ij}{\dot x}^i{\dot x}^j + e\sum_i A_i(x)
{\dot x}^i + V(x),$$ then the momentum and velocity are no
longer proportional; instead, 
\begin{equation}
\label{eq-pmagnetic}
p_i=\sum_j m_{ij}{\dot x}^j + A_i(x).
\end{equation}  

To see why (\ref{eq-pmagnetic}) and not the
simpler relation (\ref{eq-psimple}) is the appropriate definition of
momentum here, let us consider the example in
3-dimensional configuration space where
$m_{ij}=m\delta_{ij}$, the vector potential is $e x^1 {\dot x}^3$, and
$V=0$.  This lagrangian is invariant under translations of
configuration space in the direction of the $x^3$-axis; if we look for
conserved quantities, we find that
$p_3=m{\dot x}^3+ex^1$ is a constant of motion, while
the single term $m{\dot x}^3$ is not.  (The motions are circular helices
with axis in the $x^2$ direction.)

Investigating this example further, one finds that $p_2=m{\dot x}^2$ is
also a conserved quantity, while $p_1=m{\dot x}^1$ is not.  In fact,
translations along the $x^1$ axis are not symmetries of the lagrangian;
however, one may check directly from the equations of motion that there
is a third conserved quantity, namely $m{\dot x}^1-e x^3$, found
most easily from the hamiltonian formalism, as follows.
The hamiltonian $H=L-\sum_i p_i {\dot x}^i$ may be written in terms of
positions and velocities as $\frac{m}{2}\sum_i({\dot x}^i)^2$ or in
terms of positions and conjugate momenta as
$\frac{1}{2m}(p_1^2+p_2^2+(\frac{p_3}{m}-\frac{e}{m}x^1)^2). $  
   From this formula, one sees that $\frac{\partial}{\partial
x^1}+e\frac{\partial}{\partial p_3}$ is a symmetry; its generating
hamiltonian is $-ex^3+p^1=m{\dot x}^1-ex^3.$

Symplectic geometry gives a way to recover the proportionality
of momentum and velocity, with the additional gain of ``gauge
invariance,'' (i.e., the magnetic field and not the potential enters
the equations).   One 
introduces the coordinates $\xi_i=m{\dot x}^i,$ so that the
hamiltonian has the simple form $\frac{1}{2m}\sum_i(\xi_i)^2,$ but the
symplectic form $\sum dx^i\wedge dp_i$ becomes in these coordinates
$\sum_i dx^i\wedge d\xi_i + e\, dx^3\wedge dx^1.$  The $x^i$ and $\xi_i$
are no longer canonically conjugate variables; the extra term $e\,
dx^3\wedge dx^1 $ in the symplectic form is the magnetic field.  In
this formulation, the symmetry of the system with respect to rotation around
the $x^2$ axis becomes manifest.

In the language of the late 20th century, one says that the phase
space of a particle in a magnetic field is no longer the cotangent
bundle of configuration space; rather it is a ``shifted cotangent
bundle'' obtained by beginning with  the cotangent bundle of a principal
circle bundle over configuration space and then applying the operation
of symplectic reduction.  (See \cite{bo-br:differentialII} for an
exposition.)  

Another link between this example and 20th century mathematics comes
through the work of Noether \cite{no:invariante}, who, in the course
of developing ideas of Einstein and Klein in general relativity
theory, found a very general equivalence between symmetries and
conservation laws in field theory (i.e. variational problems which may
involve several independent variables) now known as {\em Noether's
theorem}.  Noether's work on this subject is so important that her
name has become inextricably attached to the relation between
symmetry and conservation laws.  It would be interesting to trace the
early history of this relation, perhaps even as far as ancient Greek
astronomy, where the uniform circular motion of heavenly bodies might
have been related to the perfection (i.e. rotational symmetry) of the
celestial sphere.

In the remainder of this paper, we will concentrate on the hamiltonian
point of view toward ``Noether's theorem.''  

\section{Hamiltonian systems with symmetry}

The relation between conserved quantities and symmetries becomes very
simple when expressed in terms of Poisson brackets.  On a symplectic
manifold $P$, the time evolution of a function $F$ under the
hamiltonian flow of a function $H$ is given by 
the
formula
$$\frac{dF}{dt}=\{F,H\}.$$  Antisymmetry of this bracket implies
immediately that $F$ is a conserved quantity for the hamiltonian flow
of $H$ if and only if the hamiltonian flow of $F$ consists of
symmetries of $H$.  If $F$ and $G$ are both conserved quantities,
then so is $\{F,G\}$, by Poisson's theorem, which
is most conveniently 
proved via the Jacobi identity.
If $F_1,\ldots,F_k$ are conserved
quantities, then so is any function $g(F_1,\ldots,F_k)$, so the
process of building new conserved quantities from old ones naturally
terminates if we arrive at a list  $F_1,\ldots,F_k$ of functions such
that, for each pair of functions in the list,
$\{F_i,F_j\}=\pi_{ij}(F_1,\ldots F_k)$ for functions $\pi_{ij}$.  Lie
\cite{li:theorie} refers to such a list of functions as generating a {\em
function group}, the function group itself consisting of all the  
$g(F_1,\ldots,F_k)$.  If $F_1,\ldots,F_k$ are functionally
independent, then the $\pi_{ij}$ are uniquely determined smooth functions
on $\reals^k$  and define on  the set of all smooth functions on
$\reals^k$ the structure of what Lie \cite{li:theorie}
 calls an {\em abstract function
group}.  This is a Lie algebra structure defined by the bracket
$$\{F,G\}=\sum_{i,j}\pi_{ij}(\mu_1,\ldots,\mu_k)\frac{\partial
F}{\partial \mu_i} \frac{\partial G}{\partial \mu_j}.$$  An
important special
case is that where the functions $\pi_{ij}$ are
linear.  The linear functions on $\reals^k$ are then closed
under the bracket operation and form a Lie algebra $\frakg$; the
carrier $\reals^k$ of the abstract function group is identified
with the dual space $\frakg^*$.

In geometric terms, a
phase space with a symmetry group consists of a manifold $P$ equipped
with a symplectic structure $\omega$ and a hamiltonian action of a Lie
group $G$.  By the latter, we mean a symplectic action of $G$ on $P$
together with an equivariant map $J$ from $P$ to the dual $\frakg^*$
of the Lie algebra of $G$ such that, {\em for each $v\in \frakg$, the
1-parameter group of transformations of $P$ generated by $v$ is the
flow of the hamiltonian vector field with hamiltonian $x\mapsto
\langle J(x),v\rangle.$}  The map $J$ is called the {\em momentum map}
(or, by many authors, {\em moment map}) of the hamiltonian action.  
  If one is simply given a
symplectic action of $G$ on $P$, any map $J$ satisfying the
condition in italics above, even if it is not equivariant, is called {\em
a} momentum map for the action.

The antisymmetry of the Poisson bracket now implies the hamiltonian
version of Noether's theorem: {\em the
momentum map $J$ is invariant under the hamiltonian flow of $H$ if and
only if $H$ is invariant under the action of the component of the
identity in the group $G$.} 
We defer until later (see Sections \ref{sec-universal} and
\ref{sec-baer})
 the question of how nonidentity components of
a symmetry group $G$ are related to conservation laws.  

The duality between symmetries and conservation laws in the
hamiltonian formalism becomes most striking in the language of {\em Poisson
geometry}, which is the contemporary name for the geometry of Lie's
abstract function groups.  We recall that a {\em Poisson manifold} is a
manifold $M$ with a Lie algebra structure $\{~,~\}$ on its algebra of
smooth functions which is an algebra derivation in each argument.
For a function $H$, the derivation $F\mapsto \{F,H\}$ is the {\em
hamiltonian vector field} associated to $H$.  If a group $G$ acts by
automorphisms of a 
symplectic manifold $P$ in such a way that the quotient space $P/G$
inherits a manifold structure, then $P/G$ becomes a Poisson manifold
through the identification of the smooth functions on $P/G$ with the
$G$-invariant smooth functions on $P$.\footnote{If $P/G$ is not 
a manifold, we may treat it as
some kind of ``Poisson variety'' or ``Poisson space.''}
  The natural projection
$\pi:P\to P/G$ is a {\em Poisson map} in the sense that its pullback
on functions is a Lie algebra homomorphism.  If the action of $G$ is
hamiltonian, its equivariant momentum map $J:P\to \frakg^*$ is also a
Poisson map, and the diagram of Poisson maps $\frakg^*\leftarrow P
\rightarrow P/G$ is called a {\em dual pair}.  In symplectic terms, this
means that the tangent spaces of the fibres of the two maps are
symplectically orthogonal to one another.  The algebras  
of functions on $\frakg^*$ and $P/G$ commute 
with one another as subalgebras of the Lie algebra of functions on
$P$, and if the fibres of the two maps in the dual pair are connected,
these two subalgebras are the full commutants of one another.

As we have suggested above, many of the ideas in this section can
already be found in Lie's book \cite{li:theorie}.  (See
\cite{we:sophus} for a brief historical discussion.)

\section{Symplectic groupoid actions}
\label{sec-groupoid}
The ``duality'' between $\frakg^*$ and $P/G$
would be more symmetrical if the fibres of $J$ could also be
seen as orbits of an action.  In general, $P/G$ is not the dual of a
Lie algebra, but there is an object which plays the role of the second
group in this picture; it is a {\em symplectic
groupoid} for the Poisson manifold $P/G$.  This is a
groupoid\footnote{We recall that a {\em groupoid} is a small category
whose morphisms are all invertible; if this category has just one
object, the groupoid is a group.} whose
objects are the points of $P/G$ and whose morphisms form a symplectic
manifold $\Gamma(P/G)$.  The symplectic structure on $\Gamma(P/G)$
is compatible with the Poisson structure on $P/G$ in the sense that
target and source maps $\alpha$ and $\beta$
from $\Gamma(P/G)$ to $P/G$ are Poisson and antiPoisson
(i.e. reversing the sign of Poisson brackets) respectively.
The groupoid multiplication is compatible with the
symplectic structure in the following sense.   Writing $\Gamma$ for
$\Gamma(P/G)$ and denoting by $\Gamma^{(2)}\subset \Gamma\times\Gamma$
the manifold
of composable pairs, we have three natural maps from $\Gamma^{(2)}$ to
$\Gamma$, namely the cartesian product projections $p_1$ and $p_2$ and
the groupoid multiplication $m$.  The symplectic structure $\omega$ on
$\Gamma$ should then satisfy the cocycle condition
$p_1^*(\omega)-m^*(\omega)+ p_2^*(\omega) = 0$. 

  The quotient 
projection $\pi:P \to P/G$ may be seen as a momentum map (in an
extended sense) for a symplectic groupoid action of $\Gamma(P/G)$ on $P$ whose
orbits are the fibres of $J$.  This means (see \cite{mi-we:moments}
that we have a
map
$a:\Gamma\times_{P/G} P\to P$ satisfying the usual
``associativity property'' of an action; the action being symplectic
means
that 
$p_1^*(\omega)-a^*(\omega_P)+p_2^*(\omega_P)=0.$   The fibre product
here is taken with respect to the maps $\beta$ and $\pi$, $\omega_P$ is
the symplectic structure on $P$, and $p_1$ and $p_2$ are again the
cartesian projections from $\Gamma\times_{P/G} P$ to $\Gamma$ and
$P$.  If we write $gx$ for $a(g,x)$, then we have the equivariance 
relation $J(gx)=gJ(x)$ involving
the natural action of the groupoid on its objects.  We call $\pi$ the
{\em momentum map} of the groupoid action.

When $\frakg$ is the Lie algebra of $G$, a natural choice for
 $\Gamma(\frakg^*)$ is
the cotangent bundle $T^*G$ with its canonical symplectic structure
and a groupoid structure whose source and target maps to
$\frakg^*=T_e^*G$ are given by left and right translation.  Under the
identification of $T^*G$ with $\frakg^* \times G$ by left
translations, the groupoid multiplication is given by $(\mu,g)(\nu,h)=\pi
(\nu,gh)$.  The fact that the fibres of $\pi:P\to P/G$ are the orbits of the
$G$ action follows from the fact that the cotangent bundle projection
is a groupoid homomorphism from $\Gamma(\frakg^*)=T^*G$ to $G$
(considered as a groupoid with a single object), and that the
symplectic groupoid action on $P$ factors through this homomorphism.
For a more general Poisson manifold like $P/G$, there is no such
homomorphism from $\Gamma(P/G)$ to a group, and so the symplectic
groupoid is absolutely essential to the picture.

Passing from groups to groupoids leads to a new definition of a
dual pair.  We may define a dual pair to
consist of a pair of symplectic groupoids $\Gamma_1$ and $\Gamma_2$
with underlying Poisson manifolds $Q_1$ and $Q_2$, together with
commuting symplectic groupoid actions of $\Gamma_1$ and $\Gamma_2$ on a
symplectic manifold $P$ having momentum maps $J_1$ and $J_2$ such that
the orbits of each action are the fibres of the momentum map of the
other.  We note that this formulation includes the case
where all the groupoids and spaces involved are discrete.
The groupoid action contains both aspects of the
symmetry/conservation duality, since the conserved quantity is just
the momentum map of the groupoid action.

Our main point, then, is that a  hamiltonian action of a  Lie group
 should be seen as a special case of 
a more flexible notion of symmetry--the action of a
symplectic groupoid.  Such an action always comes equipped with a
momentum map which is a Poisson map and which in fact determines the
action of the component of the groupoid containing the identity elements.
Pursuing this idea, it is natural to ask whether properties of
momentum maps of hamiltonian group actions extend to the groupoid
case.  As a first step toward answering this question, we have
established in \cite{we:discrete}
an extension of Kirwan's nonabelian version of Atiyah,
Guillemin, and Sternberg's convexity theorem.

\section{Poisson Lie groups and beyond}
We have seen that the cotangent bundle of a Lie group $G$ is a rather
special symplectic groupoid in that it admits a homomorphism to a
group.  This homomorphism, the cotangent bundle projection $\pi$, is a
covering morphism in the sense of \cite{br-da-ha:topological}; 
this means that its fibres
are {\em bisections} of the groupoid, i.e. they project
diffeomorphically to the space of objects under the source and target maps.
Finally, $\pi$ is a Poisson map when $G$
carries the zero Poisson structure.  This is equivalent to the fact
that the fibres of $\pi$ are lagrangian submanifolds of the cotangent
bundle.  

The properties of the symplectic groupoids $T^*G$ listed in the
paragraph above nearly characterize these examples, the only
exceptions being the quotients of cotangent bundles by lattices of
bi-invariant (and hence closed) 1-forms.  On the other hand, if we
drop the requirement that the fibres of a covering morphism $\pi$ from
a symplectic groupoid $\Gamma$ to a group $G$ be lagrangian, we obtain
a much larger class of symplectic groupoids which are still closely
connected to groups.  The group $G$ must now carry a nonzero Poisson
structure in order for $\pi$ to be a Poisson map; compatibility of the
group structure with the Poisson structure is expressed by the
condition that the multiplication $G\times G\to G$ be a Poisson map.
A symplectic groupoid action of $\Gamma$ on $P$ still corresponds to
an action of the group $G$ on $P$.  The manifold of objects of
$\Gamma$, which is the target of the momentum maps of such actions,
turns out to be a group itself; it is called the dual of the
Poisson-Lie group $G$ and is denoted $G^*$.  (The group structure on
$G^*$ is essentially encoded in the Poisson structure on $G$, and
vice-versa.)  The momentum map is equivariant with respect to an
action of $G$ on $G^*$ known as the {\em dressing action}.  Since the
fibres of $\pi$ are not lagrangian, the elements of $G$ no longer act
on $P$ by symplectic transformations.  Rather, the action as a whole
is a {\em Poisson action} in the sense that the map $G\times P\to P$
is a Poisson map.  (The dressing action is also a Poisson action.)
The original groupoid $\Gamma$ is simultaneously a symplectic groupoid
for both $G$ and $G^*$; as a symplectic manifold, it is known as the
{\em Heisenberg double} of the pair $(G,G^*)$ of dual Poisson Lie
groups.  What we have described here is essentially Lu's momentum map
theory \cite{lu:momentum} for Poisson Lie group actions.

Having dealt with compatibility conditions for
symplectic structures on groupoids and Poisson structures on groups,
one might ask for a common framework for both.  This is
provided by the notion of {\em Poisson groupoid}
\cite{we:coisotropic}.  Since Poisson structures correspond to possibly
degenerate bivector fields, one might hope for a theory which also
includes degenerate 2-forms.  This is provided by the theory of Dirac
structures \cite{co:dirac}.  These are subbundles of a direct sum
$TP\oplus T^*P$ which are maximal isotropic for a natural symmetric
bilinear form and which are closed under a bracket discovered by
Courant \cite{co:dirac} and which has become the prototype for an object
known as a {\em Courant algebroid} \cite{li-we-xu:manin}.  

Computations on moduli spaces of flat
connections in gauge theory have led to a yet more general notion of
hamiltonian symmetry.  Here, the 2 form $\omega_P$ on phase space is
neither closed nor nondegenerate, but  these ``defects'' are
compensated for by the presence of an auxiliary structure on the
group, just as the noninvariance of the symplectic structure on a
Poisson $G$-space is compensated for by the Poisson structure on the
group.  In the theory of {\em quasi-hamiltonian $G$-spaces},
introduced in \cite{al-ma-me:lie}, the
symmetry group $G$ acting on $P$ now carries a bi-invariant inner
product and, hence, a bi-invariant 3-form $\phi$ defined by
$\phi(u,v,w)=\langle [u,v],w\rangle .$ The momentum map is now {\em
group-valued}, i.e. we have $J:P\to G$ with the property that
$d\omega_P=J^*\phi$.  As in the ordinary hamiltonian case, $\omega_P$
is $G$-invariant, and $J$ is now equivariant with respect to the adjoint
action of $G$ on itself.  The relation between the momentum map and
the infinitesimal generators of the symmetry group is now more
complicated.  For $v\in \frakg$, 
the corresponding vector
field $v_P$ on $P$ satisfies the
 condition that, 
$$i_{v_P}\omega_P = J^*\langle\frac{1}{2}(\theta+\overline{\theta}),v\rangle,$$
where $\theta$ and $\overline{\theta}$ are the left-invariant and
right-invariant
Maurer-Cartan forms on $G$. Note that, since $\omega_P$ may be
degenerate, this condition does not necessarily determine the action
uniquely for a given momentum map.


The relation between momentum and symmetry becomes even looser in 
Karshon's  theory of ``abstract momentum maps'' \cite{ka:moment}.
Here, $G$ is a torus acting on $P$, and $J$ is an equivariant map from $P$ to
$\frakg^*$.  There is no longer a 2-form on $P$.  Instead, the action
and the momentum map are related by the condition that, along the
fixed point set of any subgroup $H\subseteq G$, the image of $J$ must
lie in a subspace parallel to $\frakh^{\perp}\subseteq \frakg^*$;
i.e. the component of the momentum map given by each element of $\frakh$
is a conserved quantity for the restricted action.
Remarkably, this weak relation is sufficient to prove versions of many
results usually associated with ordinary momentum maps.

\section{Universal momentum maps}
\label{sec-universal}
There have been at least two quite different approaches to the idea of
a ``universal momentum map.''  Evans and Lu
\cite{ev-lu:variety} construct for each Poisson Lie group 
a Poisson space which is the target of a momentum map 
for every Poisson action of the given group.  Ortega
and Ratiu \cite{or-ra:optimal}, on the other hand, give preeminence to the
symmetry/conservation relationship and, after fixing
a particular action, find a space whose function algebra consists of
the quantities which are conserved for all invariant 
 hamiltonian systems.  We now describe these two
approaches.  Details omitted here will appear in
 \cite{bu-we:lectures}.  

The construction by Evens and Lu may be seen as an extension of that
of Chu \cite{ch:symplectic}.  (Also see \cite{gu-st:techniques}.)
In this earlier work, one has a Lie group $G$
acting on a manifold $P$ carrying a closed (but possibly degenerate) 
2-form $\omega$ which is preserved by the action.  For each $x$ in
$P$, we may pull back $\omega$ by the map $g\mapsto gx$ from $G$ to
$P$ to get a closed, left-invariant 2-form $K(x)$ on $G$.  Such forms
correspond to 2-cocycles with constant real coefficients for the Lie
algebra $\frakg$; i.e. $K$ is a (smooth) map from $P$ to
$Z^2(\frakg)$.  The target $\frakg^*$ of usual momentum maps is the
same as the space $C^1(\frakg)$ of Lie algebra cochains, and the
coboundary operator is a natural map $\delta: \frakg^*\to
Z^2(\frakg)$.  It turns out that $Z^2(\frakg)$ is the dual of a Lie
algebra, and that $\delta$ is a Poisson map.  
If $J:P\to \frakg^*$ is an equivariant
 momentum map in the usual sense, then
the composition $\delta\circ J$ is equal to $K$. 
Notice that the map $K$ exists and is canonically associated to every
symplectic action, while $J$ may not exist and is defined only up to
an additive constant (which is in the kernel of $\delta)$.  
On the
other hand, $K$ may be much ``cruder'' than $J$.
For instance, when $\frakg$ is abelian,
$\delta$ is the zero map, and so $K$ is constant and equal to zero (i.e. the
orbits of a hamiltonian action are isotropic), even though $J$ may be
nontrivial.  

If a momentum map $J:P\to \frakg^*$ is not equivariant, then
$\delta\circ J$ differs from $K$ by a constant element $b$ of
$Z^2(\frakg)$.  If one adds the cocycle $b$ to the Poisson structure
on $\frakg^*$ as a constant term, then $J$ becomes a Poisson map and
is equivariant with respect to an affine action of $G$ on $\frakg^*$
whose linear part is the coadjoint action.  The translational part of
this action is a 1-cocycle on $G$ with values in $\frakg^*$; the
derivative of this cocycle at the identity of $G$ may be identified
with $b$. (See \cite{so:structure} or \cite{so:structureenglish}.)
One may view $\frakg^*$ with this affine Poisson structure as a
Poisson submanifold of the dual of a central extension of $\frakg$.
Although one could use the 1-dimensional extension associated to $b$,
it is interesting to use the universal central extension $X(g)$ of $\frakg$
by $Z^2(\frakg)^*$ in which the bracket $[x,y]$ of two elements of
$\frakg$ is the linear functional on $Z^2(\frakg)$ given by
$b\mapsto b(x,y)$.  There is a natural Poisson map $\hat{\delta}$
from $X(\frakg)^*$ to $Z^2(\frakg)$.  (It is not just the restriction map.)
Any momentum map $J$ for a
symplectic action of $G$ on $P$ yields an equivariant Poisson map
$\hat{J}:P\to X(\frakg)^*$ with $\hat{\delta}\circ \hat{J}=K$.
Although $\hat{J}$ is not completely determined by the action, it
does provide, for any action admitting a momentum map, a Poisson map
which lifts $K$ and which is as ``refined'' as $J$.

The definition of the map $K$ is extended to Poisson actions of a Poisson Lie group $G$
in the construction of Evens and Lu 
\cite{ev-lu:variety}.  Their target
is the space $\call(D\frakg)$ of {\em lagrangian subalgebras} of the
{\em Drinfeld double} $D\frakg$ of the Lie bialgebra $\frakg$, the Lie
algebra of the double group $DG$.  The Drinfeld
double is a Lie algebra structure on the vector space
$\frakg\oplus\frakg^*$ together with the invariant inner product
$((x,\xi),(y,\eta))\mapsto \xi(y)\eta(x);$ a lagrangian subalgebra is
a subalgebra which is maximal isotropic with respect to the inner
product.  (This is another example of a Dirac structure in a Courant
algebroid.)  Evens and Lu show that $\call(D\frakg)$ has the structure of
an algebraic variety which is stratified by Poisson manifolds
such that there is a dense subalgebra of continuous functions
on $\call(D\frakg)$ whose restrictions to the strata are smooth and
which is closed under the stratum-wise Poisson bracket.  Given a
Poisson action of $G$ on $P$, they construct a natural
equivariant map
 $EL:P\to\call(D\frakg)$ which is 
a smooth Poisson map into a single stratum.  Using the
theory of Dirac structures
\cite{co:dirac}\cite{li-we-xu:manin}\cite{li-we-xu:dirac}, it is
possible to construct the equivariant map $EL$ for arbitrary Poisson actions,
though it is generally not continuous (though it is in fact smooth
when $P$ is symplectic).  Nevertheless, $EL$ can be shown
to be a Poisson map in a certain sense. 
When $G$ has the zero Poisson structure,
the double $D\frakg$ is the semidirect product of $\frakg$ with
$\frakg^*$ via the coadjoint representation.  Among the maximal
isotropic subspaces of $\frakg\oplus\frakg^*$ are the graphs of
antisymmetric mappings $\frakg\to\frakg^*$, which correspond to skew
symmetric bilinear forms, or 2-cochains with real coefficients, on
$\frakg$.  It turns out that the subspace is a subalgebra if
and only if the cochain is a cocycle.  (The corresponding subgroup is
the graph of the $\frakg^*$-valued cocycle on $G$ connected with the
affine action discussed above.)
Thus, the target space
$Z^2(\frakg)$ is contained in $\call(D\frakg)$ (as an open
submanifold).  The Poisson structure on $Z^2(\frakg)$ induced from
the Evens-Lu structure on $\call(D\frakg)$ turns out to be identical
to that which occurs in the construction of Chu, 
and the maps $EL$ and $K$ coincide when $P$ is symplectic.

Continuing with the case where $G$ has the zero Poisson
structure, we also find among the lagrangian subalgebras of $D\frakg$
the sums $\frakh\oplus\frakh^\perp$, where $\frakh$ is an arbitrary
subalgebra of $\frakg$ and $\frakh^\perp$ is its annihilator in
$\frakg^*$.  These occur as the values of $EL$ when the manifold $P$
carries the zero Poisson structure: $EL(p)=\frakg_p\oplus
\frakg_p^\perp,$ where $\frakg_p$ is the Lie algebra of the isotropy
subgroup $G_p$ of $p$ in $G$.  (The corresponding subgroup of the
double group is the conormal bundle of $G_p$.)  Thus, in the absence
of any Poisson structure at all, we still have a ``momentum map''
which essentially assigns to each point in a $G$-manifold its isotropy
submanifold.  The level sets of this momentum map, or {\em isotropy types}, are
invariant under any $G$-equivariant map.  

Of course, there are no interesting $G$-invariant hamiltonian systems
when $P$ has the zero Poisson structure, since all hamiltonian vector
fields are zero; nevertheless, this example provides a nice lead-in
to the construction of Ortega and Ratiu \cite{or-ra:optimal}, since
their construction is motivated in part by the fact that conventional
momentum maps may miss the ``conservation of isotropy''.  

Ortega and Ratiu begin with a group $G$ acting by automorphisms on a
Poisson manifold $M$.  Let $\cale$ be the set of hamiltonian vector
fields of $G$-invariant functions defined on $G$-invariant open
subsets of $M$.  The (local) flows of the vector fields in $\cale$
form a pseudogroup $\calg$ of Poisson diffeomorphisms between open
subsets of $M$, and we may form the {\em momentum space} $M/\calg$.
In general, $M/\calg$ is not a manifold, but it is a ``smooth Poisson
space'' in the sense that, if we define the smooth functions on
$M/\calg$ to be the $\calg$-invariant smooth functions on $M$, the
space of these functions forms a Poisson algebra.  Ortega and Ratiu's
{\em optimal momentum map} is simply the natural projection map
$\calj:M\to M/\calg$.  Immediately from its definition, it has the
property of being being conserved by the flow of any $G$-invariant
hamiltonian.  Since the elements of $\calg$ commute with $G$, $G$ acts
in a natural way on the momentum space, and the optimal momentum map
is equivariant.  It is universal in the sense that, given any
equivariant Poisson map $K:M\to P$ which is conserved by the flows of
all $G$-invariant hamiltonians, there is a unique equivariant Poisson
map $\phi:M/\calg\to P$ such that $K=\phi\circ\calj$.  The optimal
momentum map can be nontrivial even if a momentum map in the usual
sense does not exist.  For hamiltonian actions, it may be more refined
than the usual momentum map $J$.  For instance, in the case of a
proper hamiltonian action, each level set of $\calj$ is a connected
component of the intersection of a level set of $J$ with an isotropy
type.

We remark that the optimal momentum map is very closely related to the
  Poisson $\Gamma$ structure introduced by Condevaux, Dazord, and
   Molino \cite{co-da-mo:geometrie} in a paper whose title is
   essentially the same as that of this one.  

\section{The Baer groupoid}
\label{sec-baer}
We can forge closer links among some of the concepts in this paper with
the aid of the {\em Baer groupoid} $\calb G$ associated to a group
$G$.  The set of objects of this groupoid is the set $\cals G$
of subgroups of $G$; a
morphism from $H$ to $K$ is a subset of $G$ which is at the same time
a right coset $Ka$ and a left coset $bH$.  For such a morphism, it
follows that  $bH=aH$, so that $K=aHa^-1$, and $K$ and $H$ are
conjugate.  The product of two morphisms is the usual product of
subsets in a group.

Thus, the $\calb G$ orbits in $\cals G$ are
the conjugacy classes, and the isotropy subgroup
$H\in \cals G$ is isomorphic to the quotient $NH/H$, where $NH$ is
the normalizer of $H$ in $G$.  \footnote{When $G$ is a compact
connected Lie group, the restriction of $\calb G$ to the set of
maximal tori is a transitive groupoid which might be called the {\em
Weyl groupoid}, since its isotropy groups are the Weyl groups.  A
possible advantage of this groupoid over the individual Weyl groups is
that it admits a nontrivial action of the automorphism group of $G$.}

Any action of a group $G$ on a space $P$ induces an action of $\calb
G$ on $P$.  The momentum map of this action assigns to each $x\in P$ the
isotropy subgroup $G_x$, and the coset $aG_x$ maps $x$ to $ax$.  The
level sets of this momentum map are exactly the isotropy types, so we
are close to the Evens--Lu momentum map.  If the
action is hamiltonian, we can extend the $\calb G$ action on $P$ to an
action of the product groupoid $\calb G \times T^*G$ with object space
$\cals G \times \frakg^*$.  The momentum map of this extended action has
as level sets the intersection of the level sets of $J:P\to \frakg^*$
with the isotropy types.   The connected components of these level
sets are, for proper actions, precisely the level sets of the optimal
momentum map.  Except for this question of connectedness, we have
recovered the optimal momentum map using a natural construction
involving a momentum space which depends only on the group 
and not on the particular (proper) action under consideration.

Further development of these ideas looks like a good project for
the early 21st century.

\end{document}